\documentclass[a4paper,12pt,reqno]{amsart}

\newtheorem{thm}{Theorem}[section]

\newtheorem{prop}[thm]{Proposition}

\theoremstyle{definition}
\newtheorem{remark}[thm]{Remark}

\numberwithin{equation}{section}

\def\al{\alpha}
\def\be{\beta}
\def\ga{\gamma}
\def\de{\delta}
\def\ep{\varepsilon}

\def\la{\lambda}

\def\si{\sigma}
\def\vp{\varphi}

\def\La{\Lambda}
\def\Ga{\Gamma}

\def\Z{\mathbb{Z}}

\def\R{\mathbb{R}}
\def\C{\mathbb{C}}
\def\N{\mathbb{N}}

\def\cD{\mathcal D}

\def\cP{\mathcal P}

\def\lc{\text{\rm lc}}

\newcommand{\rphis}[5]{\,_{#1}\vp_{#2} \left( \genfrac{.}{.}{0pt}{}{#3}{#4}
\ ;#5 \right)}

\newcommand{\rFs}[5]{\,_{#1}F_{#2} \left( \genfrac{.}{.}{0pt}{}{#3}{#4}
\ ;#5 \right)}

\title[Spectral properties]{Spectral properties of operators \\ using tridiagonalisation}
\author{Mourad E.H. Ismail and Erik Koelink}
\date{\today}
\address{City University of Hong Kong, Department of Mathematics, 83 Tat Chee Avenue, Kowloon, Hong Kong
and King Saud University, Riyadh, Saudi
Arabia.}
\email{meismail@cityu.edu.hk}
\address{Radboud Universiteit Nijmegen, IMAPP, FNWI, Heyendaalseweg 135, 6525 AJ Nijmegen, the Netherlands}
\email{e.koelink@math.ru.nl}
 
\begin{document}
\date{\today}
\begin{abstract} A general scheme for tridiagonalising differential, difference or $q$-difference operators using orthogonal polynomials is described. From the tridiagonal form the spectral decomposition can be described in terms of the orthogonality measure of generally different orthogonal polynomials. 
Three examples are worked out: (1) related to Jacobi and Wilson polynomials for a second order differential operator, (2) related to little $q$-Jacobi polynomials and Askey-Wilson polynomials for a bounded second order $q$-difference operator, (3) related to little $q$-Jacobi polynomials for an unbounded second order $q$-difference operator. In case (1) a link with the Jacobi function transform is established, for which we give a $q$-analogue using example (2). 
\end{abstract}

\maketitle

%%%%%%%%%%%%%%%%%%%%%%%%%%%%%%%%%%%%%%%%%%%%%%%%%%%%%%%%%%%%%%%%%%%%%%%%%%%%%%%%%%%%%%
%%%%%%%%%%%%%%%%%%%%%%%%%%%%%%%%%%%%%%%%%%%%%%%%%%%%%%%%%%%%%%%%%%%%%%%%%%%%%%%%%%%%%%
%%%%%%%%%%%%%%%%%%%%%%%%%%%%%%%%%%%%%%%%%%%%%%%%%%%%%%%%%%%%%%%%%%%%%%%%%%%%%%%%%%%%%%

\section{Introduction}

Bochner's 1929 theorem classifies polynomials which are eigenfunctions to a second order differential operator, see e.g. \cite{Isma} where also a historical discussion and extensions to other operators can be found. Instead of looking for differential operators diagonalised by polynomials, we can also look for differential operators which are tridiagonalised by polynomials. Naturally, there is no need to restrict to differential operators, and we can also consider difference and $q$-difference operators. Neither is there a need to restrict to a  differential or difference operator of second-order. In \cite{IsmaK} we have discussed a general framework as well as the examples of the Schr\"odinger operator with Morse potential and the Lam\'e equation. In this paper we give a general approach for finding self-adjoint differential (or difference) operators which are tridiagonal in a suitable basis of orthogonal polynomials.  Since the spectral decomposition of tridiagonal operators can be expressed completely in terms of, in general different, orthogonal polynomials, we obtain the spectral decomposition of the differential or difference operator involved. 

In Section \ref{sec:genframework} we sketch the general framework, and we specialise to the case where the basis in which the operator is tridiagonal is given by orthogonal polynomials. Since we start with orthogonal polynomials which are eigenfunctions to an  explicit differential or difference operator, it is natural to start from the families in the ($q$-)Askey scheme. 
In Section \ref{sec:exJacobi} we work out the case of the Jacobi polynomials, which yields an operator which is, after a change of variable, for a special choice the differential operator for the Jacobi functions. Its spectrum can be given explicitly since the tridiagonal operator corresponds to Wilson polynomials. In Section \ref{sec:exlittleqJacobi} we discuss a $q$-analogue of this situation involving the little $q$-Jacobi polynomials. However, due to lack of symmetry in the $q$-case we have to consider two cases. It turns out that one case can be solved within the $q$-Askey scheme using Askey-Wilson polynomials, whereas the other case leads to orthogonal polynomials not in the $q$-Askey scheme. 
In Section \ref{sec:exJacobi} we link to a result \cite[(3.3)]{Koor-LNM} 
by Koornwinder which states that Jacobi polynomials can be mapped to Wilson polynomials by the Jacobi function transform, and we give a $q$-analogue of this result in Section \ref{ssec:lqJacobi-case1}. Koornwinder's result is predated by a result by Broad \cite{Broa}, which the limit case of the Whittaker transform mapping Laguerre polynomials to Meixner-Pollaczek polynomials, and is hugely extended by Groenevelt \cite{Groe} to the top of the Askey-scheme by showing that the Wilson transform maps Wilson polynomials to Wilson polynomials. 

Starting within the Askey-scheme and requiring that the tridiagonalisation can be solved within the Askey-scheme leaves only the cases of the Jacobi and Laguerre polynomials. As stated, the first case is discussed in Section  \ref{sec:exJacobi}, whereas the second is a similar calculation which can also be obtained using limit transitions from the results in Section \ref{sec:exJacobi}. The other cases in the Askey-scheme to which this procedure can be applied are the Meixner, continuous Hahn, continuous dual Hahn and Wilson polynomials. The orthogonal polynomials associated to the tridiagonalisation are no longer in the Askey-scheme, but these are orthogonal polynomials associated to birth and death processes which is a general phenomenon. A further study of these examples is desirable, and may lead to a better understanding of this procedure. A similar remark applies to the orthogonal polynomials in the $q$-Askey scheme, see also Section \ref{sec:conclrmk}.

The tridiagonalisation method is also known as the $J$-matrix method. 
Tridiagonalisation for different operators in all kinds of bases is a much employed method in many physical problems, see e.g. \cite{AlHaBAAA}, \cite{Hayd}, \cite{YamaR-PhysRevA}, and references in \cite{IsmaK}. We hope that the more general approach given here widens the class of operators to which the tridiagonalisation method applies. In Section \ref{sec:conclrmk} we discuss briefly the possible extension of this method to an even larger class of operators. 

%%%%%%%%%%%%%%%%%%%%%%%%%%%%%%%%%%%%%%%%%%%%%%%%%%%%%%%%%%%%%%%%%%%%%%%%%%%%%%%%%%%%%%
%%%%%%%%%%%%%%%%%%%%%%%%%%%%%%%%%%%%%%%%%%%%%%%%%%%%%%%%%%%%%%%%%%%%%%%%%%%%%%%%%%%%%%
%%%%%%%%%%%%%%%%%%%%%%%%%%%%%%%%%%%%%%%%%%%%%%%%%%%%%%%%%%%%%%%%%%%%%%%%%%%%%%%%%%%%%%

\section{General framework}\label{sec:genframework}

%%%%%%%%%%%%%%%%%%%%%%%%%%%%%%%%%%%%%%%%%%%%%%%%%%%%%%%%%%%%%%%%%%%%%%%%%%%%%%%%%%%%%%
%%%%%%%%%%%%%%%%%%%%%%%%%%%%%%%%%%%%%%%%%%%%%%%%%%%%%%%%%%%%%%%%%%%%%%%%%%%%%%%%%%%%%%
%%%%%%%%%%%%%%%%%%%%%%%%%%%%%%%%%%%%%%%%%%%%%%%%%%%%%%%%%%%%%%%%%%%%%%%%%%%%%%%%%%%%%%

In this section we introduce the general set-up for tridiagonalising certain differential or difference operators. In Section \ref{ssec:genframework} we give the general framework which is specialised to the case of orthogonal polynomials in Section \ref{ssec:orthopols}. 

\subsection{General framework}\label{ssec:genframework}

Let $\mu$ and $\nu$ be positive Borel measures  on the real line $\R$
so that $\nu$ is absolutely continuous with respect to $\mu$. 
Let $\de = \frac{d\nu}{d\mu}$ be the Radon-Nikodym derivative, so $\de\geq 0$.
Let $\cD$ be a function space that is dense in both $L^2(\mu)$ and $L^2(\nu)$
with respect to the respective topologies. 

We assume that we have an operator 
$L \colon \cD\to \cD$ operator such that 
$(L,\cD)$ is symmetric in $L^2(\mu)$. Moreover, we assume the 
existence of an orthonormal basis $(\Phi_n)_{n\in\N}$ of $L^2(\mu)$
of eigenfunctions of $L$, so $\Phi_n\in \cD$ for all $n\in\N$ and
$L\Phi_n = \La_n \Phi_n$, $\La_n\in\R$. Then $(L,\cD)$ is essentially self-adjoint.
% So the closure of $(L,\cD)$ is
% a self-adjoint operator $(\bar{L},D(\bar{L}))$ with 
% \begin{equation}\label{eq:defmaxdom}
% D(\bar{L}) \, = \, \{ f = \sum_{n=0}^\infty f_n\, \Phi_n\in L^2(\mu)\, \, \colon \, \, \sum_{n=0}^\infty 
% |f_n \La_n|^2 < \infty \, \}, \quad 
% f_n = \langle f, \Phi_n\rangle_{L^2(\mu)}. 
% \end{equation}

Next we assume that there is an orthonormal basis $(\phi_n)_{n=0}^\infty$ 
of $L^2(\nu)$ such that $\phi_n\in \cD$ for all $n\in\N$ and that
\begin{equation}\label{eq:phicombinationofPhi}
\phi_n \, = \, A_n\, \Phi_n \, + \, B_n \, \Phi_{n-1}, \qquad A_n, B_n\in \R
\end{equation}
(with the convention $B_0=0$). Finally we assume that 
multiplication by the reciprocal of the Radon-Nikodym
derivative preserves the space $\cD$, so $M\colon \cD \to \cD$, 
$Mf=\de^{-1}f$.  

The object of study is the operator $T=ML\colon \cD\to \cD$ as 
an operator on $L^2(\nu)$.
Note 
\begin{equation}\label{eq:actionT0onphin}
\begin{split}
\langle T\phi_n, \phi_m\rangle_{L^2(\nu)} \, =\, & 
\langle A_n\, ML \Phi_n \, + \, B_n \, ML\Phi_{n-1}, 
A_m\, \Phi_m \, + \, B_m \, \Phi_{m-1} \rangle_{L^2(\nu)} \\
 =\, & \langle A_n\, \La_n \Phi_n \, + \, B_n \, \La_{n-1} \Phi_{n-1}, 
A_m\, \Phi_m \, + \, B_m \, \Phi_{m-1} \rangle_{L^2(\mu)} \\
 =\, & \La_n A_n  B_{n+1}\de_{n+1,m} \, + \, (A_n^2\, \La_n +B_n^2 \La_{n-1}) \de_{n,m}
\, + \, \La_{n-1} A_{n-1}  B_{n}\de_{n,m+1}
\end{split}
\end{equation}
In particular, $(T, \cD)$ is a Jacobi operator on the 
Hilbert space $L^2(\nu)$; 
\begin{equation}\label{eq:Tastridiagop}
\begin{split}
&T\phi_n \, = \, a_n\, \phi_{n+1}\, + \, b_n\, \phi_n \, + \, a_{n-1}\, \phi_{n-1}, \\
a_n\, &=\, \La_n\, A_n \, B_{n+1}, \qquad b_n\, = \, \La_n \, A_n^2 \, + \, \La_{n-1} \, B_n^2,
\end{split}
\end{equation}
where the convention $B_{0}=0$ gives $a_{-1}=0$.  Since $a_n, b_n\in \R$ and  
$(\phi_n)_{n=0}^\infty$  is an orthonormal basis of $L^2(\nu)$ it follows
that $T$ with domain $D(T)$ finite linear combinations of 
the orthonormal basis $(\phi_n)_{n=0}^\infty$  is a symmetric densely defined
operator to which standard techniques of tridiagonal (or Jacobi) operators
can be applied, see e.g. \cite{Akhi}, \cite{Bere}, \cite{DunfS}, \cite{IsmaK}, \cite{Koel}, \cite{Simo} and references given there.
We assume $D(T) = \cD$ so that we can obtain the spectral decomposition of suitable self-adjoint extensions
of $(T, \cD)$. We assume that $(T,\cD)$ is essentially self-adjoint. 
In general the spectral decomposition of the closure of $(T,\cD)$  depends on the possible vanishing of the 
coefficients $a_n$ and on the growth behaviour of the sequences 
$(a_n)_{n=0}^\infty$, $(b_n)_{n=0}^\infty$, see e.g. \cite[Thm.~2.7]{IsmaK}. 

\begin{remark}\label{rmk:Laddingconstant} Note that we can trivially change $L$ by adding a suitable real constant, say $L^\ga = L+\ga$, and we change $\La_n$ to $\La_n^\ga = \La_n + \ga$. The operator 
in general changes by adding a multiplication operator $T^\ga = T + \ga M$. By comparing the coefficients $a_n^\ga=a_n + \ga A_nB_{n+1}$, $b_n^\ga=b_n + \ga (A_n^2+B_n^2)$ we find the three-term recurrence 
\[
M\, \phi_n \, = \, A_nB_{n+1} \phi_{n+1} \, + \, (A_n^2+B_n^2)\phi_n \, + \, 
A_{n-1}B_{n} \phi_{n-1}.
\]
The additional freedom is useful, see Sections \ref{sec:exJacobi} and \ref{sec:exlittleqJacobi}. 
\end{remark}

%%%%%%%%%%%%%%%%%%%%%%%%%%%%%%%%%%%%%%%%%%%%%%%%%%%%%%%%%%%%%%%%%%%%%%%%%%%%%%%%%%%%%%
%%%%%%%%%%%%%%%%%%%%%%%%%%%%%%%%%%%%%%%%%%%%%%%%%%%%%%%%%%%%%%%%%%%%%%%%%%%%%%%%%%%%%%
%%%%%%%%%%%%%%%%%%%%%%%%%%%%%%%%%%%%%%%%%%%%%%%%%%%%%%%%%%%%%%%%%%%%%%%%%%%%%%%%%%%%%%

\subsection{Specialising to orthogonal polynomials}\label{ssec:orthopols}

The examples we consider are based on orthogonal polynomials. 
Let $\mu$ and $\nu$ be orthogonality measures for orthogonal 
polynomials;
\begin{equation}\label{eq:munuorthomeasures}
\int_\R P_n(x)P_m(x)\, d\mu(x)\,=\, H_n\de_{n,m}, \qquad
\int_\R p_n(x)p_m(x)\, d\nu(x)\,=\, h_n\de_{n,m}.
\end{equation}
We assume that both $\mu$ and $\nu$ correspond to a determinate 
moment problem, so that the space $\cP$ of polynomials is dense in $L^2(\mu)$ and $L^2(\nu)$. 
% ??could also be $N$-extremal; examples known??.
We also assume that $\int_\R f(x)\, d\mu(x) = \int_\R f(x)r(x)\, d\nu(x)$, 
where $r$ is a polynomial of degree $1$, 
so that the Radon-Nikodym derivative 
$\frac{d\nu}{d\mu}=\de = 1/r$. 
Then we obtain, using $\lc(p)$ for the leading coefficient of a polynomial $p$, 
\begin{equation}\label{eq:explexpgenpols}
p_n\, = \, \frac{\lc(p_n)}{\lc(P_n)}\, P_n \, + \, \lc(r)
\frac{h_n}{H_{n-1}}\frac{\lc(P_{n-1})}{\lc(p_n)}\, P_{n-1}
\end{equation}
and by taking $\phi_n$, respectively $\Phi_n$, the corresponding 
orthonormal polynomials to $p_n$, respectively $P_n$, we see that
\eqref{eq:phicombinationofPhi} specialises to 
\begin{equation}\label{eq:OPphicombinationofPhi}
\begin{split}
&\phi_n \, = \, A_n\, \Phi_n \, + \, B_n \, \Phi_{n-1}, \qquad A_n \, = \, \frac{\lc(p_n)}{\lc(P_n)} \sqrt{\frac{H_n}{h_n}}, \quad
B_n \, = \, \lc(r)
\sqrt{\frac{h_n}{H_{n-1}}}\frac{\lc(P_{n-1})}{\lc(p_n)}.
\end{split}
\end{equation}
We assume the existence of a symmetric self-adjoint operator $L$ with domain $\cD=\cP$ on $L^2(\mu)$ with $LP_n =\La_n P_n$, and so $L\Phi_n =\La_n \Phi_n$, for eigenvalues $\La_n\in \R$. Then $(L,\cP)$ is an essentially self-adjoint operator, since $\cP$ is dense in $L^2(\mu)$, as assumed in Section \ref{ssec:genframework}. By convention $\La_{-1} =0$. So this means that $(P_n)_n$ satisfies a bispectrality property, and we can take the family $(P_n)_n$ from the Askey scheme or its $q$-analogue. 

Considering the operator $T = rL$ with domain $\cD = \cP$ on $L^2(\nu)$ shows that
the coefficients in \eqref{eq:Tastridiagop} are given by 
\begin{gather*}
a_n \, = \, \La_n \lc(r) \frac{\lc(p_n)}{\lc(p_{n+1})} \sqrt{\frac{h_{n+1}}{h_n}}, 
\quad
b_n \, = \, \La_n \frac{H_n}{h_n} \left( \frac{\lc(p_n)}{\lc(P_n)}\right)^2 
\, + \, \La_{n-1} \lc(r)^2 \frac{h_n}{H_{n-1}} 
\left( \frac{\lc(P_{n-1})}{\lc(p_n)}\right)^2
\end{gather*}

By switching from $L$ and $T$ to $L^\ga$ and $T^\ga$, see Remark \ref{rmk:Laddingconstant}, for a suitable constant $\ga$ we can assume that  
$\La_n^\ga \not=0$ for all $n\in \N$. 
\begin{equation}\label{eq:Tonorthocompl}
T^\ga \phi_n\, = \, a_n^\ga \phi_n \, +\,   b_n^\ga \phi_n\,  +\,  a_{n-1}^\ga \phi_{n-1},
\end{equation}
and $a_n^\ga \not=0$ for $n\geq 0$, and $a_{-1}^\ga=0$.  
So we need to solve for the orthonormal polynomials $r_n(\la)$ 
% and corresponding monic polynomials $R_n(\la)$ 
satisfying 
\begin{equation}\label{eq:recurrenceforspectralpols}
\begin{split}
\la r_n(\la) &= a_n^\ga r_n(\la) +  b_n^\ga r_n(\la) + a_{n-1}^\ga r_{n-1}(\la), 
% \\
% \la R_n(\la) &= R_n(\la) +  \be_n R_n(\la) + \al_{n-1}^2 R_{n-1}(\la).
\end{split}
\end{equation}
We assume that the corresponding Jacobi operator is essentially self-adjoint (or equivalently that the orthogonal polynomials $\{r_n\}$ correspond to a determinate moment problem). Thus 
we have a positive Borel measure $\rho$ so that 
\begin{equation}\label{eq:orthrelpolrn}
\int_\R r_n(\la) r_m(\la)\, d\rho(\la) = \de_{n,m}.
\end{equation} 
In particular, this implies that we assume $T$ with domain $\cP$ essentially self-adjoint. 

\begin{remark}\label{rmk:bdpsupport} Polynomials associated to a birth and death process with birth-rate $\la_n$ and death-rate $\mu_n$ are generated by 
\[
-x F_n(x) = \la_{n} F_{n+1}(x) - (\la_n+\mu_n) F_n(x) + \mu_n F_{n-1}(x), \quad 
\ F_0(x)=1,\ F_1(x) = (\la_0+\mu_0 -x)/\la_0. 
\]
with $\la_n>0$ for all $n$ and $\mu_n>0$ for $n\geq 1$, $\mu_0\geq 0$, see \cite[Ch. 5]{Isma}. 
Apart from sign issues the  recursion \eqref{eq:Tonorthocompl} corresponds to a birth and death process with birth-rate $\la_n = |\La_{n}^\ga| \frac{H_{n}}{h_{n}} \left( \frac{\lc(p_{n})} {\lc(P_{n})}\right)^2$ and death-rate 
$\mu_n = |\La_{n-1}^\ga| (\lc(r))^2 \frac{h_{n}}{H_{n-1}} \left( \frac{\lc(P_{n-1})}{\lc(p_{n})}\right)^2$ for $n>0$, $\mu_0=0$. Assuming that the signs of $a_n^\ga$ and $\La_n^\ga$ are independent of $n\in\N$, this allows us to draw conclusions on the support of the spectrum of $\bar{T^\ga}$, i.e. the orthogonality measure of the polynomials $r_n$, see \cite[Ch. 5, 7]{Isma}. In this case, let $\ep=\text{sgn}(a_n^\ga)$, $\eta=\text{sgn}(\La_n^\ga)$. 
In particular, the spectrum of $\bar{T}$ is contained in $[0,\infty)$ if
$(\ep,\eta)=(-1,1)$ or $(1,1)$ and the spectrum of $\bar{T}$ is contained in $(-\infty,0])$ if $(\ep,\eta)=(1,-1)$ or $(-1,-1)$. In the examples $\La_n^\ga$ is bounded from above or from below, so that it is possible to find a $\ga\in \R$ so that all eigenvalues $\La_n^\ga$ have the same sign. In the examples we see that the condition on $a_n^\ga$ is also valid. 
\end{remark}

Defining $U\colon L^2(\nu) \to L^2(\rho)$, $U\colon \phi_n \mapsto r_n$, gives a unitary operator satisfying $UT = U^\ast M$, where $M$ is the multiplication (by $\la$) operator on $L^2(\rho)$. It follows that the spectrum of $T$ equals $\si(T) = \text{supp}(\rho)$, and the spectrum is simple. 

It follows that 
\begin{equation}\label{eq:formaleigvector}
F(x;\la) = \sum_{n=0}^\infty r_n(\la) \phi_n(x) =
\sum_{n=0}^\infty r_n(\la) \frac{p_{n}(x)}{\sqrt{h_{n}}} 
\end{equation}
is a formal eigenvector for the eigenvalue $\la$ for $T$. The series in \eqref{eq:formaleigvector} converges in $L^2(\nu)$ for $\la$ a discrete mass point of $\rho$ and diverges for $\la$ in the continuous part of $\rho$. In that case the series, which is a non-symmetric Poisson kernel,  
\begin{equation}\label{eq:genPoissonkernel}
F_t(x;\la) = 
\sum_{n=0}^\infty t^n r_n(\la) \frac{p_{n+1}(x)}{\sqrt{h_{n+1}}} 
\end{equation}
converges for $|t|<1$.

%%%%%%%%%%%%%%%%%%%%%%%%%%%%%%%%%%%%%%%%%%%%%%%%%%%%%%%%%%%%%%%%%%%%%
%%%%%%%%%%%%%%%%%%%%%%%%%%%%%%%%%%%%%%%%%%%%%%%%%%%%%%%%%%%%%%%%%%%%%
%%%%%%%%%%%%%%%%%%%%%%%%%%%%%%%%%%%%%%%%%%%%%%%%%%%%%%%%%%%%%%%%%%%%%
%%%%%%%%%%%%%%%%%%%%%%%%%%%%%%%%%%%%%%%%%%%%%%%%%%%%%%%%%%%%%%%%%%%%%

\section{Jacobi polynomials}\label{sec:exJacobi}

We work out the programme of Section \ref{ssec:orthopols} for the case of the Jacobi polynomials and the related hypergeometric differential operator. 

For the Jacobi polynomials $P^{(\al,\be)}_n(x)$, we follow the standard notation 
 \cite{AndrAR}, \cite{Isma}, 
\cite{KoekS}. %, \cite{Rain}. 
We take the measures $\mu$ and $\nu$ to be the
orthogonality measures for the Jacobi polynomials for parameters $(\al+1,\be)$, 
and $(\al,\be)$ respectively. We assume $\al,\be>-1$. So we set $P_n(x)=P^{(\al+1,\be)}_n(x)$, 
$p_n(x)=P^{(\al,\be)}_n(x)$. This gives
\[
\begin{split}
h_n &= N_n(\al) = \frac{2^{\al+\be+1}}{2n+\al+\be+1}
\frac{\Ga(n+\al+1)\Ga(n+\be+1)}{\Ga(n+\al+\be+1)\, n!}, \quad H_n = N_n(\al+1), 
\\ 
\lc(p_n) &= l_n(\al) 
= \frac{(n+\al+\be+1)_n}{2^n n!}, \quad \lc(P_n) = l_n(\al+1).
\end{split}
\]
Moreover, $r(x)=1-x$. 
Note that we could have also shifted in $\be$, but due to the symmetry 
$P_n^{(\al,\be)}(x) = (-1)^n P_n^{(\be,\al)}(-x)$ 
of the Jacobi polynomials in $\al$ and $\be$ it suffices to consider the
shift in $\al$ only.

The Jacobi polynomials are eigenfunctions of 
\begin{equation*}
\begin{split}
&L f (x)\, = \, (1-x^2)\, f''(x) + \bigl( \be-\al-1-(\al+\be+3)x\bigr) f'(x), \\
&LP^{(\al+1,\be)}_n \, = \, -n(n+\al+\be+2)\, P^{(\al+1,\be)}_n
\end{split}
\end{equation*}
so that $\La_n = -n(n+\al+\be+2)$. We set $\ga= -(\al+\de+1)(\be-\de+1)$, so that
$\La_n^\ga = -(n+\al+\de+1)(n+\be -\de+1)$.
Hence,  we study on $L^2([-1,1], (1-x)^\al (1+x)^\be\, dx)$ the operator 
$(1-x)(L+\ga)$ which is
\begin{equation}\label{eq:JacobidefT}
T^\ga\,  = \, (1-x)(1-x^2)\, \frac{d^2}{dx^2} + (1-x)\bigl( \be-\al-1-(\al+\be+3)x\bigr) \frac{d}{dx} - (1-x)(\al+\de+1)(\be-\de+1).
\end{equation}

Applying the results of Section \ref{ssec:orthopols} gives the explicit expression for the recursion coefficients
\begin{equation*} %\label{eq:Jacobianbn}
\begin{split}
a_n^\ga &=  \frac{2 (n+\al+\de+1)(n+\be-\de+1)}{2n+\al+\be+2} 
\sqrt{\frac{(n+1)\, (n+\al+1)\, (n+\be+1)\, (n+\al+\be+1)}
{(2n+\al+\be+1)\, (2n+\al+\be+3)}}
\\ 
b_n^\ga &=  - \frac{2(n+\al+\de+1)(n+\be-\de+1)(n+\al+1)(n+\al+\be+1)}
{(2n+\al+\be+1)\, (2n+\al+\be+2)} \\ & \qquad
\, - \, \frac{2 n (n+\be) (n+\al+\de+1)(n+\be-\de)}
{(2n+\al+\be)\, (2n+\al+\be+1)}. 
\end{split}
\end{equation*}
Then the recursion relation \eqref{eq:Tonorthocompl} for $\frac12 T^\ga$ is solved by the orthonormal version of the Wilson polynomials
\[
W_n(\mu^2; \frac12(1+\al), \frac12(1+\al)+\de, \frac12 (1-\al)+\be-\de, \frac12 (1+\al)),
\]
see \cite{KoekS}, where the relation between the eigenvalue $\la$ of $T$ and $\mu^2$ is given by $\la = -2 \left(\frac{\al+1}{2}\right)^2-2\mu^2$. 

\begin{thm}\label{thm:Jacobi} Let $\al>-1$, $\be>-1$, and assume $\ga= -(\al+\de+1)(\be-\de+1)\in \R$. The unbounded operator $(T^\ga, \cP)$ defined by \eqref{eq:JacobidefT} on  $L^2([-1,1], (1-x)^\al (1+x)^\be\, dx)$
with domain the polynomials $\cP$  is essentially self-adjoint. The spectrum  
of the closure $\bar{T^\ga}$ is simple and given by 
\[
\begin{split}
(-\infty, &-\frac12(\al+1)^2) \cup 
\{ -\frac12(\al+1)^2+2(\frac12(1+\al)+\de+k)^2 \, \colon \,  k\in\N, \ \frac12(1+\al)+\de+k<0\} \\ & \cup 
\{ -\frac12(\al+1)^2+2(\frac12(1-\al)+\be-\de+l)^2 \, \colon \,  l\in\N, \ \frac12(1-\al)+\be-\de+l<0\} 
\end{split}
\]
where the first set gives the absolutely continuous spectrum and the other sets correspond to the discrete spectrum of the closure of $T^\ga$. The discrete spectrum is either empty or at most one of these sets is non-empty. 
\end{thm}

\begin{remark}\label{rmk:thmJacobi} Note that in Theorem \ref{thm:Jacobi} $\de\in \R$ or $\Re \de = \frac12(\be-\al)$. In the second case there is no discrete spectrum. 
\end{remark}

The eigenvalue equation $T^\ga f_\la = \la f_\la$ is a second-order differential operator with regular singularities at $-1$, $1$, $\infty$. In the Riemann-Papperitz notation it 
is 
\[
\mathcal{P} \left\{ \begin{matrix} -1 & 1 & \infty &  \\
0 & -\frac12(1+\al) + i \tilde\la & \al+\de+1 & x \\
-\be & -\frac12(1+\al) + i \tilde\la & \be-\de+1 &  
                    \end{matrix}
\right\}
\]
with the reparametrisation $\la = -\frac12(\al+1)^2 -2\tilde\la^2$ of the spectral parameter. In case $\ga=0$, i.e. $\de=-\al-1$ or $\de=1+\be$,
  We set $\de=1+\be$, and we put $T=T^0$, $a_n^0=a_n$, $b_n^0=0$. Note that we have degenerate case, since $a_0=0$ as well as $b_0=0$ and $a_n\not=0$ for $n\geq 1$. It follows that $T\phi_0= b_0\phi_0 =0$, and we have to study the recurrence relation \eqref{eq:Tonorthocompl} on $(\C\phi_0)^\perp$, i.e. we start at $n=1$ instead of at $n=0$. The recurrence relation \eqref{eq:Tonorthocompl} starting from $n=1$ can also be solved similarly in terms of 
Wilson polynomials, namely by the orthonormal version of
\[
W_n(\mu^2; \frac12(1+\al), \frac12(1-\al), \frac12 (3+\al), \be + \frac12(3+\al)),
\]
where the relation between the eigenvalue $\la$ of $T$ and $\mu^2$ is given by $\la = -2 \left(\frac{\al+1}{2}\right)^2-2\mu^2$ as in the derivation of Theorem \ref{thm:Jacobi}. 

\begin{prop}\label{prop:Jacobi} The operator $(T,\cP)$ defined by \eqref{eq:JacobidefT} with $\ga=0$  with domain the polynomials
$\cP\subset L^2([-1,1], (1-x)^\al (1+x)^\be\, dx)$ is essentially self-adjoint. The spectrum  
of the closure $\bar{T}$ is simple and given by 
\[
(-\infty, -\frac12(\al+1)^2) \cup 
\{ -\frac12(\al+1)^2+2(\frac12(1-\al)+k)^2 \, \colon \,  k\in\N, \ \frac12(1-\al)+k<0\} \cup \{0\}
\]
where the first set gives the absolutely continuous spectrum and the other sets correspond to the discrete spectrum of the closure of $T$. 
\end{prop}

\begin{remark}\label{rmk:propJacobi} Note that putting $\de=1+\be$ (or $\de=-1-\al$) in the result of Theorem \ref{thm:Jacobi} precisely gives back the statement of Proposition \ref{prop:Jacobi}. The spectral decomposition $U$ is now defined on $(\C\phi_0)^\perp$ by $\phi_{n+1}\mapsto r_n$, where the $r_n$ are the orthonormal Wilson polynomials with parameters $\frac12(1+\al), \frac12(1-\al), \frac12 (3+\al), \be + \frac12(3+\al))$. 
\end{remark}

The eigenvalue equation $T f_\la = \la f_\la$ is a second-order differential operator with regular singularities at $-1$, $1$ and $\infty$, 
hence can be solved in terms of hypergeometric functions. Changing variables 
$x=1-2 \cosh^{-2}t$ transforms $2T$ in the second-order differential operator
\[
\frac{d^2}{dt^2} + \bigl( (2\be+1)\coth t \, - \, (2\al+2\be +3)\tanh t\bigr) \frac{d}{dt}
\]
on $L^2((0,\infty), 2^{\al+\be+2} (\sinh t)^{2\be+1} (\cosh t)^{-2\al-2\be-3}\, dt)$. This is the Jacobi differential operator as in \cite{Koor-JF}, where $(\al,\be)$ of     \cite{Koor-JF} correspond to $(\be, -\al-\be-2)$. Proposition \ref{prop:Jacobi} provides an alternative proof of the simplicity and location of the spectrum of the Jacobi differential operator, see \cite[Thms. 2.3, 2.4]{Koor-JF}, but it does not give the Jacobi functions as eigenfunctions. 

Having established the link of $T$ to the Jacobi differential operator, we know that the Jacobi function transform gives the spectral decomposition. In particular, applying the integral transformation 
\[
 \hat{f}(\la) = \int_0^\infty \phi_\la^{(\be,-\al-\be-2)} f(t) (2\sinh t)^{2\be+1} (2\cosh t)^{-2\al-2\be-3}\, dt 
\]
for $f\in L^2((0,\infty),  (2\sinh t)^{2\be+1} (2\cosh t)^{-2\al-2\be-3}\, dt)$ see \cite{Koor-JF}, \cite{Koor-LNM} for details, to \eqref{eq:Tonorthocompl} starting from $n=1$ shows that the Jacobi function transform of $t\mapsto \phi_{n+1}(1-\cosh^{-2}t)$ satisfies the same recurrence relation as the orthonormal Wilson polynomials corresponding to Remark \ref{rmk:propJacobi}. Since the solution to the recurrence relation is completely determined by its starting value we
obtain 
\begin{gather*}
\int_0^\infty \phi^{(\be, -\al-\be-2)}_\la (t) P^{(\al,\be)}_{n+1}(1-2\cosh^{-2}t)\,
 (2\sinh t)^{2\be+1} (2\cosh t)^{-2\al-2\be-3}\, dt \, = \, \\
\sqrt{\frac{h_{n+1}S_0}{h_1 S_n}} W_n(\frac{\la^2}{4}; \frac12(1+\al), \frac12(1-\al), \frac12 (3+\al), \be + \frac12(3+\al)) \\
\times \int_0^\infty \phi^{(\be, -\al-\be-2)}_\la (t) P^{(\al,\be)}_{1}(1-2\cosh^{-2}t)\,
 (2\sinh t)^{2\be+1} (2\cosh t)^{-2\al-2\be-3}\, dt,
\end{gather*}
where $S_n$ is the squared norm of the Wilson polynomials. However, this is a special case of a result by Koornwinder \cite[\S 9]{Koor-JF}, \cite[(3.3)]{Koor-LNM} mapping Jacobi polynomials to Wilson polynomials. In our case there is a drop in the degree of the polynomial, which does not occur in \cite[(3.3)]{Koor-LNM}. Applying \cite[(3.3)]{Koor-LNM} to the integral on the left hand side gives a Wilson polynomial of degree $n+1$ with parameters $(\frac12(\al+1),\frac12(\al+1), -\frac12(\al+1),
\be +\frac12(\al+3))$. Since two of the parameters add up to zero the first term of the hypergeometric series vanishes and the polynomial can be written as a linear term times a  Wilson polynomial of degree $n$. This can also be viewed as a special case of the so-called `Diophantine' properties introduced by Calogero et al., see \cite[(3.10)]{ChenI} and references given there. 

\begin{remark}\label{rmk:JacobilimitWhittaker} The Whittaker (confluent hypergeometric) differential operator can be obtained as a limit case of the Jacobi differential operator, see \cite{Koor-LNM}. The operator arises from this approach if we start with Laguerre polynomials instead of Jacobi polynomials.
\end{remark}

%%%%%%%%%%%%%%%%%%%%%%%%%%%%%%%%%%%%%%%%%%%%%%%%%%%%%%%%%%%%%%%%%%%%%
%%%%%%%%%%%%%%%%%%%%%%%%%%%%%%%%%%%%%%%%%%%%%%%%%%%%%%%%%%%%%%%%%%%%%
%%%%%%%%%%%%%%%%%%%%%%%%%%%%%%%%%%%%%%%%%%%%%%%%%%%%%%%%%%%%%%%%%%%%%
%%%%%%%%%%%%%%%%%%%%%%%%%%%%%%%%%%%%%%%%%%%%%%%%%%%%%%%%%%%%%%%%%%%%%

\section{Little $q$-Jacobi polynomials}\label{sec:exlittleqJacobi}

Next we work out the programme of Section \ref{ssec:orthopols} for the case of the little $q$-Jacobi polynomials as a $q$-analogue of Section \ref{sec:exJacobi}. 

We follow the notation of \cite{KoekS}. For $0<a<q^{-1}$, $b<q^{-1}$ we define the little $q$-Jacobi polynomials by 
\[
p_n(x;a,b;q) = \rphis{2}{1}{q^{-n}, abq^{n+1}}{aq}{q,qx}.
\]
The little $q$-Jacobi polynomials are orthogonal with respect to a discrete measure on $q^\N$;
\begin{equation}\label{eq:ortholittleqJacobipols}
\begin{split}
&\sum_{k=0}^\infty w_k(a,b;q)  p_n(q^k;a,b;q) p_m(q^k;a,b;q) = \de_{n,m} h_n(a,b;q), \\
&w_k(a,b;q) = \frac{(qb;q)_k}{(q;q)_k} (aq)^k, \quad
h_n(a,b;q) = \frac{(abq^2;q)_\infty}{(aq;q)_\infty} \frac{(1-abq) \, (aq)^n}{(1-abq^{2n+1})} 
\frac{(q,qb;q)_n}{(qa,qab;q)_n}
\end{split}
\end{equation}
The leading coefficient is given by
\[
\text{lc}\left( p_n(\cdot;a,b;q)\right) = (-1)^n q^{-\frac12 n(n-1)} \frac{(abq^{n+1};q)_n}{(aq;q)_n}.
\]

Defining the second-order $q$-difference operator
\begin{gather*}
(L^{(a,b)}f)(x) = \frac{B(x)}{x} \left( f(qx)-f(x)\right) + 
\frac{D(x)}{x} \left( f(q^{-1}x)-f(x)\right), 
\\ B(x) = a(bqx-1), \qquad D(x) = x-1
\end{gather*}
we have 
\[
L^{(a,b)} p_n(\cdot;a,b;q)\, = \, \La_n(a,b;q)\, p_n(\cdot;a,b;q), \qquad \La_n(a,b;q) = q^{-n} (1-q^n)(1-abq^{n+1}).
\]

%%%%%%%%%%%%%%%%%%%%%%%%%%%%%%%%%%%%%%%%%%%%%%%%%%%%%%%%%%%%%%%%%%%%%
%%%%%%%%%%%%%%%%%%%%%%%%%%%%%%%%%%%%%%%%%%%%%%%%%%%%%%%%%%%%%%%%%%%%%
%%%%%%%%%%%%%%%%%%%%%%%%%%%%%%%%%%%%%%%%%%%%%%%%%%%%%%%%%%%%%%%%%%%%%
%%%%%%%%%%%%%%%%%%%%%%%%%%%%%%%%%%%%%%%%%%%%%%%%%%%%%%%%%%%%%%%%%%%%%
\subsection{Case one: shift in $a$}\label{ssec:lqJacobi-case1}

In the notation of Section \ref{ssec:orthopols} we take
$P_n(x) = p_n(x;aq,b;q)$, $H_n = h_n(aq,b;q)$, $p_n(x) = p_n(x;a,b;q)$, 
$h_n=h_n(a,b;q)$. We assume $0<a<q^{-1}$ and $b<q^{-1}$. Then $\mu$ and $\nu$ are discrete measures and
$\mu(\{q^k\})= w_k(aq,b;q) = q^k w_k(a,b;q) = q^k \nu(\{q^k\})$, so that 
$r(x) = x$. Also $L= L^{(aq,b)}$, $\La_n = \La_n(aq,b;q)$ and $T^\ga$ is given as
a second order difference operator on $L^2(\nu)$ by 
\begin{equation}\label{eq:lqJp1-T}
(T^\ga f)(x) = %\bigl(x(L+\ga)f\bigr)(x) =  
aq(bqx-1) \left( f(qx)-f(x)\right) + 
(x-1) \left( f(q^{-1}x)-f(x)\right) + \ga x f(x). 
\end{equation}
In particular, $T^\ga$ is a bounded operator on $L^2(\nu)$. 

We now replace $\ga = (1+abq^2) -(acq+bq/c)$, so that 
\[
\La_n^\ga = q^{-n} (1-acq^{n+1})(1-bq^{n+1}/c),
\]
and we assume $c$ choosen so that $\ga\in\R$ and $\La_n^\ga>0$
% \in\R\setminus\{0\}$ 
for all $n\in \N$. This is the case when e.g. $bq<c< 1/aq$, and by the conditions on $a$, $b$ this is always possible. 
A calculation gives the values
\begin{equation}\label{eq:lqJp1-anbn}
\begin{split}
a_n^\ga &= -\frac{\sqrt{aq} (1-acq^{n+1})(1-bq^{n+1}/c)}{(1-abq^{2n+2})}
\sqrt{\frac{(1-q^{n+1})(1-aq^{n+1})(1-bq^{n+1})(1-abq^{n+1})}{(1-abq^{2n+1})(1-abq^{2n+3})}}, \\
b_n^\ga &= \frac{(1-acq^{n+1})(1-aq^{n+1})(1-abq^{n+1})(1-bq^{n+1}/c)}{(1-abq^{2n+1})(1-abq^{2n+2})}
 \\ &\qquad\qquad 
+ \frac{aq(1-q^n)(1-acq^{n})(1-bq^{n})(1-bq^{n}/c)}{(1-abq^{2n})(1-abq^{2n+1})}.
\end{split}
\end{equation}

The corresponding recurrence relation \eqref{eq:recurrenceforspectralpols} for 
$\frac{-1}{\sqrt{aq}}T^\ga$ corresponds exactly to the orthonormal Askey-Wilson polynomials with parameters $(\sqrt{aq}, c\sqrt{aq}, \frac{b}{c}\sqrt{\frac{q}{a}},\sqrt{aq})$. Explicitly,
\[
\begin{split}
r_n(\la) &= C \, p_n(x;\sqrt{aq}, c\sqrt{aq}, \frac{b}{c}\sqrt{\frac{q}{a}},\sqrt{aq}|q) \\
C & = (acq,\frac{bq}{c};q)_n \sqrt{\frac{(1-abq) ( q,aq,bq;q)_n}{(1-abq^{2n+1}}}
\end{split}
\]
where $\la = 1+aq -2\sqrt{aq} x$ and using the standard notation for the Askey-Wilson polynomials \cite{AndrAR}, \cite{GaspR}, \cite{Isma}, \cite{KoekS}. So the continuous spectrum of $T$ is the interval $[1+aq-2\sqrt{aq}, 1+aq+2\sqrt{aq}]=[(1-\sqrt{aq})^2,(1+\sqrt{aq})^2]$. 
Discrete spectrum in the orthogonality measure can only occur if $|c\sqrt{aq}|>1$ or $|\frac{b}{c}\sqrt{\frac{q}{a}}|>1$, see \cite{GaspR}, \cite{KoekS}. Using the explicit form of the discrete mass points in the orthogonality measure for the Askey-Wilson polynomials then leads to the following theorem. 

\begin{thm}\label{thm:lqJacobi1-bdd} Consider the operator $T^\ga$ defined by \eqref{eq:lqJp1-T} with $\ga = (1+abq^2) -(acq+bq/c)\in\R$ and with the assumption $\La_n^\ga>0$  for all $n\in\N$. Then $T^\ga$ acting on
$L^2(\nu)$ with $\nu$ the discrete orthogonality measure  \eqref{eq:ortholittleqJacobipols} for the little $q$-Jacobi polynomials with $0<qa<1$, $b<q^{-1}$ is a bounded self-adjoint operator with simple spectrum at 
\[
\begin{split}
[(1-\sqrt{aq})^2,(1+\sqrt{aq})^2] 
&\cup \{(1-q^{-k}/c)(1-acq^{1+k}) \, \colon \,  |q^{k}c\sqrt{aq}| >1, \ k\in\N\} \\
&\cup  \{ (1-bq^{l+1}/c)(1-acq^{-l}/b)\, \colon \,   |q^l \frac{b}{c}\sqrt{\frac{q}{a}}|>1, \ l\in \N\}
\end{split}
\]
\end{thm}

\begin{remark}\label{rmk:thmlqJacobi1-bdd} Theorem \ref{thm:lqJacobi1-bdd} is a $q$-analogue of Theorem \ref{thm:Jacobi}. However in Theorem \ref{thm:lqJacobi1-bdd} we can have that the discrete spectrum can be empty, consist of one of these sets or consists of both sets. The degenerate case $\ga=0$ or $c=1/aq$ can be considered as a $q$-analogue of the Jacobi differential operator, cf. Proposition \ref{prop:Jacobi}. 
\end{remark}

We can also obtain the spectral decomposition of $T$ acting on $L^2(\nu)$ in a direct way. 
Suppose $y_\la\colon q^\N\to \C$ satisfies $T^\ga y_\la= \la y_\la$. Put
$y_\la (q^k) = \frac{p_k(\la)}{\sqrt{w_k(a,b;q)}}$,
then being an eigenfunction of $T$ is equivalent to 
\begin{gather*}
\left( \frac{aq+1-\la}{\sqrt{aq}}\right) p_k(\la) = 
\sqrt{(1-q^{k+1})(1-bq^{k+1})}p_{k+1}(\la) \\
\qquad\qquad +q^k(c\sqrt{aq} +\frac{b}{c}\sqrt{\frac{q}{a}}) p_k(\la) +
\sqrt{(1-q^{k})(1-bq^{k})}p_{k-1}(\la).
\end{gather*}
This corresponds to the orthonormal Al-Salam--Chihara polynomials with $(c\sqrt{aq}, \frac{b}{c}\sqrt{\frac{q}{a}})$ as parameters, so 
\begin{equation}\label{eq:pklaisAlSalamChihara}
p_k(\la) =  \frac{1}{\sqrt{(q,bq;q;)_k}}Q_k(x; c\sqrt{aq}, \frac{b}{c}\sqrt{\frac{q}{a}}|q)
\end{equation}
where $2x = (aq+1-\la)/\sqrt{aq}$. Here we use the notation for the Al-Salam--Chihara polynomials \cite{GaspR}, \cite{Isma}, \cite{KoekS}. The spectrum equals the support of the orthogonality measure $\si$ of the Al-Salam--Chihara polynomials $Q_k(x; c\sqrt{aq}, \frac{b}{c}\sqrt{\frac{q}{a}}|q)$, and we obtain another proof of Theorem 
\ref{thm:lqJacobi1-bdd}. 

So the spectral decomposition of $T^\ga$ on $L^2(\nu)$ is given by the unitary map 
$V \colon L^2(\nu) \to L^2(\si)$
\begin{equation}\label{eq:lqJ1spectraldecomp}
\begin{split}
(Vf)(x) &= \langle f, y_\la\rangle_{L^2(\nu)} = \sum_{k=0}^\infty f(q^k) \, \frac{Q_k(x;c\sqrt{aq}, \frac{b}{c}\sqrt{\frac{q}{a}}|q)}{\sqrt{(q,bq;q)_k}}\, \sqrt{w_k(a,b;q)} \\
&= \sum_{k=0}^\infty f(q^k) \, \frac{(\sqrt{aq})^k}{(q;q)_k }Q_k(x;c\sqrt{aq}, \frac{b}{c}\sqrt{\frac{q}{a}}|q)  
\end{split}
\end{equation} 
where $\la = aq+1-2\sqrt{aq}x$, $x=\cos \theta$. The unitarity of $f$ can also be checked directly using the orthogonality relations for the Al-Salam--Chihara polynomials. 
So $VT^\ga =MV$, where $V$ is the multiplication by $aq+1-2\sqrt{aq}\,x$. 

Applying $V$ to the three-term recurrence for $T^\ga$ gives 
\[
\bigl( aq+1-2\sqrt{aq}\,x\bigr) (V\phi_{n})(x) = a_n (V\phi_{n+1})(x) 
+ b_n (V\phi_{n})(x) + a_{n-1} (V\phi_{n-1})(x),
\]
so that $V\phi_n$ is a constant multiple of the orthonormal Askey-Wilson polynomials with parameters $\sqrt{aq}, c\sqrt{aq}, \frac{b}{c}\sqrt{\frac{q}{a}},\sqrt{aq}$. The constant is $V\phi_0$. Explicitly, 
\begin{equation}\label{eq:lqJ1-transform}
\begin{split}
&(Vp_{n}(\cdot;a,b;q))(x) = 
\sum_{k=0}^\infty p_{n}(q^k;a,b;q)) \, \frac{(aq)^{k/2}}{(q;q)_k}
Q_k(x;c\sqrt{aq}, \frac{b}{c}\sqrt{\frac{q}{a}}|q) \\ = 
& \frac{(aq)^{n/2} (acq^{n+1}, bq^{n+1}/c;q)_\infty}{(aq;q)_n\, (\sqrt{aq} e^{i\theta},\sqrt{aq} e^{-i\theta};q)_\infty }
\, p_n(x; \sqrt{aq}, c\sqrt{aq}, \frac{b}{c}\sqrt{\frac{q}{a}}, \sqrt{aq}|q)
\end{split}
\end{equation}
where $x=\cos\theta$. Since $V$ is unitary 
$h_n\de_{n,m} = \langle V\bigl(p_n(\cdot;a,b;q)\bigr),V\bigl(p_m(\cdot;a,b;q)\bigr)\rangle_{L^2(\si)}$ establishes the Askey-Wilson 
polynomials $p_n(x; \sqrt{aq}, c\sqrt{aq}, \frac{b}{c}\sqrt{\frac{q}{a}}, \sqrt{aq}|q)$ as orthogonal polynomials  with respect to 
\[
\frac{d\si}{(\sqrt{aq} e^{i\theta},\sqrt{aq} e^{-i\theta},\sqrt{aq} e^{i\theta},\sqrt{aq} e^{-i\theta};q)_\infty}
\]
which is well-known \cite{AndrAR}, \cite{GaspR}, \cite{Isma}, \cite{KoekS}.

It is straightforward to prove \eqref{eq:lqJ1-transform} 
from the explicit expression for the little $q$-Jacobi polynomials and the generating function \cite{GaspR}, \cite[(15.1.10)]{Isma}, \cite[(3.8.13)]{KoekS} for the Al-Salam--Chihara polynomials, this leads to the slightly more general result 
\begin{equation}\label{eq:generaltransform}
\sum_{k=0}^\infty \frac{t^k}{(q;q)_k}p_n(q^k;a,b;q)  Q_k(x;c,d|q) 
= \frac{(tc,td;q)_\infty}{(te^{i\theta}, te^{i\theta};q)_\infty}
\rphis{4}{3}{q^{-n}, abq^{n+1}, te^{i\theta}, te^{-i\theta}}{aq, tc, td}{q,q}  
\end{equation}
for $|t|<1$ and $x=\cos\theta$. 

In order to get the full four-parameter Askey-Wilson polynomials as orthogonal polynomials from this interpretation we observe that $f_n \colon q^k \mapsto \al^k p_n(q^k;a,b;q)$ and
$g_n \colon q^k \mapsto \al^{-k} p_n(q^k;a,b;q)$ for $\sqrt{aq}<\al<1/\sqrt{aq}$ are biorthogonal funtions in $L^2(\nu)$. By \eqref{eq:generaltransform} we find
\begin{equation}\label{eq:qKoornwinder1985}
\bigl( V f_n\bigr)(x) = 
 \frac{\al^n (aq)^{n/2} (\al acq^{n+1}, \al bq^{n+1}/c;q)_\infty}{(aq;q)_n\, (\al \sqrt{aq} e^{i\theta},\al \sqrt{aq} e^{-i\theta};q)_\infty }
\, p_n(x; \al \sqrt{aq}, c\sqrt{aq}, \frac{b}{c}\sqrt{\frac{q}{a}}, \sqrt{aq}/\al|q)
\end{equation}
and a similar expression for $Vg_n$. Now $h_n\de_{nm} = \langle Vf_n, Vg_m\rangle_{L^2(\si)}$ gives the full four-parameter family of Askey-Wilson polynomials 
$p_n(x; \al \sqrt{aq}, c\sqrt{aq}, \frac{b}{c}\sqrt{\frac{q}{a}}, \sqrt{aq}/\al|q)$ as orthogonal polynomials with respect to 
\[
\frac{d\si}{(\al \sqrt{aq} e^{i\theta},\al \sqrt{aq} e^{-i\theta},\al^{-1}\sqrt{aq} e^{i\theta},\al^{-1}\sqrt{aq} e^{-i\theta};q)_\infty}.
\]
So we have given yet another proof of the well-known orthogonality relations of the Askey-Wilson polynomials, see \cite{AndrAR}, \cite{GaspR}, \cite{Isma}, \cite{KoekS}. 

\begin{remark}
By the results of Section \ref{sec:exJacobi} we may consider \eqref{eq:qKoornwinder1985} 
as a $q$-analogue of \cite[\S 9]{Koor-JF}, \cite[(3.3)]{Koor-LNM} mapping the (bi-)orthogonal little $q$-Jacobi polynomials to the full four-parameter family of Askey-Wilson polynomials, cf. Remark \ref{rmk:thmlqJacobi1-bdd}. 
To see the formal limit transition of the Al-Salam--Chihara polynomials (with $c=1/aq$ corresponding to $\ga =0$) of \eqref{eq:pklaisAlSalamChihara} to the Jacobi functions $\phi_\mu^{(\be,-\al-\be-2)}$ we replace $e^{i\theta}=q^{\frac12 i \mu}$, $\mu \in [0, 2\pi/\ln q]$. Substituting as well $a=q^\al$, $b=q^\be$, we can write the Al-Salam--Chihara polynomial of \eqref{eq:pklaisAlSalamChihara} as
\begin{gather*}
Q_k(\frac12 (q^{\frac12 i\mu}+q^{-\frac12 i\mu}); q^{-\frac12(1+\al)}, q^{\be+1+\frac12(1+\al)} |q) 
= \\ \sum_{l=0}^\infty \frac{(q^{\frac12(-1-\al-i\mu)},q^{\frac12(-1-\al+i\mu)};q)_l}
{(q, q^{\be+1};q)_l} q^l \Bigl( (z^{-1};q)_l \Bigr)\vert_{z=q^k} 
\stackrel{\longrightarrow}{q\uparrow 1} \\
\rFs{2}{1}{\frac12(-1-\al-i\mu), \frac12(-1-\al+i\mu)}{\be+1}{1-\frac{1}{z}}
= \phi_\mu^{(\be,-\al-\be-2)}(t), \qquad z= \cosh^{-2}t
\end{gather*}
This is only a formal limit transition, and we refer to \cite{Koor1990} for another limit transitions of orthogonal polynomials from the $q$-Askey scheme to Jacobi functions. Koornwinder \cite{Koor1990} gives a rigorous limit transition from the continuous $q$-ultraspherical polynomials to Jacobi functions with $\al=\be$, where the degree of the polynomials is related to the argument of the Jacobi function as in the formal limit described in this remark.  The Jacobi function transform is related to the spherical Fourier transform on rank one symmetric spaces, see \cite{Koor-JF}, and for the case $a=b=1$, and $c=1/aq$ this corresponds to the occurrence of the Al-Salam--Chihara polynomials in the spectral decomposition of the Laplace-Beltrami operator on the quantum disc, see Vaksman \cite[\S 1.5.4]{Vaks}. 
\end{remark}

%%%%%%%%%%%%%%%%%%%%%%%%%%%%%%%%%%%%%%%%%%%%%%%%%%%%%%%%%%%%%%%%%%%%%
%%%%%%%%%%%%%%%%%%%%%%%%%%%%%%%%%%%%%%%%%%%%%%%%%%%%%%%%%%%%%%%%%%%%%
%%%%%%%%%%%%%%%%%%%%%%%%%%%%%%%%%%%%%%%%%%%%%%%%%%%%%%%%%%%%%%%%%%%%%
%%%%%%%%%%%%%%%%%%%%%%%%%%%%%%%%%%%%%%%%%%%%%%%%%%%%%%%%%%%%%%%%%%%%%
\subsection{Case two: shift in $b$}\label{ssec:lqJacobi-case2}

In the notation of Section \ref{ssec:orthopols} we take
$P_n(x) = p_n(x;a,bq;q)$, $H_n = h_n(a,bq;q)$, $p_n(x) = p_n(x;a,b;q)$, 
$h_n=h_n(a,b;q)$. We assume $0<a<q^{-1}$ and $b<q^{-1}$. Then $\mu$ and $\nu$ are discrete measures and
$\mu(\{q^k\})= w_k(a,bq;q) = \frac{(1-bq^{k+1})}{(1-bq) }w_k(a,b;q) = r(q^k) \nu(\{q^k\})$, where 
$r(x) = \frac{1-bqx}{1-bq}$. Also $L= L^{(a,bq)}$, $\La_n = \La_n(a,bq;q)$ we take $\ga=(1+abq^2)-(acq+bq/c)$ as in Section \ref{ssec:lqJacobi-case1}. Again  we assume 
$\La_n^\ga>0$ for all $n\in\N$.  Then $T^\ga$ is given as
a second order difference operator on $L^2(\nu)$ by 
\begin{equation}\label{eq:lqJp2-T}
\begin{split}
(T^\ga f)(x) =  & \frac{a(bq^2x-1)(1-bqx)}{x\, (1-bq)} \left( f(qx)-f(x)\right) + \\ & 
\frac{(x-1)(1-bqx)}{x\, (1-bq)} \left( f(q^{-1}x)-f(x)\right) 
+ \big((1+abq^2)-(acq+bq/c)\bigr) \frac{(1-bqx)} {(1-bq)} f(x)
\end{split}
\end{equation}
Note that $T$ is unbounded as an operator on $L^2(\nu)$. 

Then we can calculate 
\begin{equation}
\begin{split}
a_n^\ga &= \frac{bq\sqrt{aq}}{1-bq}\frac{(1-acq^{n+1})(1-bq^{n+1}/c)}{(1-abq^{2n+2})}
\sqrt{\frac{(1-q^{n+1})(1-aq^{n+1})(1-bq^{n+1})(1-abq^{n+1})}{(1-abq^{2n+1})(1-abq^{2n+3})}}\\
b_n^\ga &= q^{-n}\frac{(1-acq^{n+1})(1-bq^{n+1}/c)(1-abq^{n+1})(1-bq^{n+1})}
{(1-bq)(1-abq^{2n+1})(1-abq^{2n+1})} \\ & \qquad\qquad\quad + 
ab^2q^{n+2} \frac{(1-acq^{n+1})(1-bq^{n+1}/c)(1-q^{n})(1-aq^{n})} {(1-bq)(1-abq^{2n})(1-abq^{2n+1})}
\end{split} 
\end{equation}
Note that $(a_n^\ga)_n$ is a bounded sequence, even though $T^\ga$ is an unbounded operator. 
Since the sequence $(b_n^\ga)_n$ is a sum of a bounded sequence and unbounded sequence tending to infinity, it follows that \cite[Thm. VII1.4]{Bere} gives that $T^\ga$ is essentially self-adjoint. However, these polynomials do not fit in the $q$-Askey-scheme. This proves the first statement of the following proposition. The statement on the location of the spectrum follows, since we are dealing with birth and death processes, see Remark \ref{rmk:bdpsupport} since we assume $\La_n^\ga>0$. 

\begin{prop}\label{prop:lqJ-2statement}
The unbounded operator $T$ defined by \eqref{eq:lqJp2-T} on
$L^2(\nu)$ with $\nu$ the discrete orthogonality measure  \eqref{eq:ortholittleqJacobipols} for the little $q$-Jacobi polynomials with $0<qa<1$, $b<q^{-1}$ is an unbounded self-adjoint operator with simple spectrum 
located in $[0,\infty)$. 
\end{prop}

As in Section \ref{ssec:lqJacobi-case1} we can try to solve the eigenvalue equation $Ty_\la = \la y_\la$ directly. Putting $y_\la(q^k) = p_k(\la)/\sqrt{w_k(a,b;q)}$ leads to 
the recurrence 
\begin{equation}\label{eq:lqJacobi2-directeigfunctions}
\begin{split}
&(1-bq)\la\, p_k(\la) = 
a_k p_{k+1}(\la) \, + \, b_k p_k(\la) \, + \, a_{k-1} p_{k-1}(\la), \\
& a_k = -\sqrt{a} q^{-k-\frac12} (1-bq^{k+2}) \sqrt{1-bq^{k+1})(1-q^{k+1})}, \\
& b_k = q^{-k}(1-bq^{k+1})\bigl( a(1-cq^{k+1})+(1-bq^{k+1}/c)\bigr).
\end{split}
\end{equation}
Using \cite[Thm. VII1.4]{Bere} we see that the polynomials correspond to a determinate moment problem. 
It is remarkable that the support of the orthogonality measure of \eqref{eq:lqJacobi2-directeigfunctions} and the one arising 
from Proposition \ref{prop:lqJ-2statement} coincide. 

%%%%%%%%%%%%%%%%%%%%%%%%%%%%%%%%%%%%%%%%%%%%%%%%%%%%%%%%%%%%%%%%%%%%%
%%%%%%%%%%%%%%%%%%%%%%%%%%%%%%%%%%%%%%%%%%%%%%%%%%%%%%%%%%%%%%%%%%%%%
%%%%%%%%%%%%%%%%%%%%%%%%%%%%%%%%%%%%%%%%%%%%%%%%%%%%%%%%%%%%%%%%%%%%%
%%%%%%%%%%%%%%%%%%%%%%%%%%%%%%%%%%%%%%%%%%%%%%%%%%%%%%%%%%%%%%%%%%%%%

\section{Concluding remarks}\label{sec:conclrmk}

We comment on several possible extensions of the results described, which are outside the scope of this paper. 

% First of all, it is desirable to have an explicit spectral decomposition of the second order differential operator \eqref{eq:JacobidefT} for $\ga\not= 0$. This might lead to a generalisation of Koornwinder's result \cite[(3.3)]{Koor-LNM}.
% % ??Can the eigenfunctions be obtained from the Frobenius method??
% 
% Secondly, the method can be worked out for other examples within the Askey-scheme for which the conditions of Section \ref{ssec:orthopols} hold. These contain at least  
% the Meixner, continuous Hahn, continuous dual Hahn, and Wilson polynomials. However, the resulting recurrence relation \eqref{eq:Tonorthocompl} cannot be matched with known polynomials in the Askey-scheme. The coefficients in \eqref{eq:Tonorthocompl} are higher order rational functions in $n$, and could be of interest from the point of view of birth and death processes, cf. Remark \ref{rmk:bdpsupport}. More examples from the $q$-Askey scheme have yet to be explored. 

In the context of Section \ref{sec:genframework} we can also obtain tridiagonalisations of higher order differential or difference operators by looking at the operator $p(L) M p(L)$ for $p$ a real-valued function, such as a polynomial of degree $d$ and $M$ a multiplication by a linear function. Then $T$ is again a symmetric differential or ($q$-)difference operator on $L^2(\mu)$ if $L$ is self-adjoint on $L^2(\mu)$ and diagonal with a suitable basis of eigenfunctions. Note that the degree of $T$ is $2d$ on $L^2(\mu)$. A possible link to Krall-type polynomials needs to be investigated, see e.g. \cite{EverL}.

Secondly, we can extend the relation \eqref{eq:phicombinationofPhi} to allow for more terms. This leads to higher-order recurrences for $T$, which may be solved by matrix-valued orthogonal polynomials. Conversely, for known operators $T$ we can obtain spectral information on matrix-valued orthogonal polynomials that occur in this way. A similar remark is valid if we look for tridiagonalisation on $\Z$ instead of $\N$. 

\medskip
\emph{Acknowledgement.} 
The research of Mourad E.H. Ismail is supported by a Research Grants Council of Hong Kong  under contract \# 101410 and a grant from King Saud
University, Saudi Arabia.

This work was partially supported by a grant from the `Collaboration Hong Kong - Joint Research Scheme' sponsored by the Netherlands Organisation of Scientific Research and the Research Grants Council fo Hong Kong (Reference number: 600.649.000.10N007).


\begin{thebibliography}{99}

\bibitem{Akhi} N.I.~Akhiezer, 
\emph{The Classical Moment Problem and Some Related Questions in Analysis},
Oliver and Boyd, 1965. 

\bibitem{AlHaBAAA} A.D.~Alhaidari, H.~Bahlouli, M.S.~Abdelmonem, F.~Al-Ameen, T. Al-Abdulaal, {Regularization in the $J$-matrix method of scattering revisited},  Phys. Lett. A  {364}  (2007) 372--377.

\bibitem{AndrAR}  G.E.~Andrews, R.A.~Askey, R.~Roy,
\emph{Special Functions}, Cambridge Univ. Press, 1999.

\bibitem{Bere}
J.M.~Berezanski\u\i,
\emph{Expansions in Eigenfunctions of Selfadjoint Operators},
Transl. Math. Monographs 17, AMS, 1968.

\bibitem{Broa} J.T.~Broad, \emph{Extraction of continuum properties from $L^2$ basis set matrix representations of the Schr\"odinger equation: the Sturm sequence polynomials and Gauss quadrature}, in: 
\emph{Numerical Integration of Differential Equations and Large Linear Systems}, (ed. J.~Hinze),
LNM 968, Springer, 1982, pp. 53--70.


\bibitem{ChenI} Y.~Chen, M.E.H.~Ismail, \emph{Hypergeometric origins of Diophantine properties associated with the Askey scheme}, Proc. Amer. Math. Soc. 138, 943--951. 

\bibitem{DunfS}
N.~Dunford and J.T.~Schwartz,
\emph{Linear Operators II: Spectral Theory},
Interscience, 1963.

\bibitem{EverL} W.N.~Everitt, L.L.~Littlejohn, \emph{Orthogonal polynomials and spectral theory: a survey}, in:  \emph{Orthogonal Polynomials and Their Applications (Erice, 1990)}, (eds. C.~Brezinski, L.~Gori, A.~Ronveaux), IMACS Ann. Comput. Appl. Math. 9, Baltzer, 1991,  pp. 21--55

\bibitem{GaspR} G.~Gasper and M.~Rahman,
\emph{Basic Hypergeometric Series}, 2nd ed.,
Cambridge Univ. Press, 2004.


\bibitem{Groe} W.~Groenevelt, \emph{The Wilson function transform}, 
Int. Math. Res. Not. {2003}, 2779--2817. 

\bibitem{Hayd} R.~Haydock, {The recursion method and the Schr\"odinger equation}, in: P.~Nevai (Ed.),
{Orthogonal Polynomials}, NATO Adv. Sci. Inst. Ser. C Math. Phys. Sci. 294, Kluwer, 1990, 
 pp. 217--228 

\bibitem{Isma} M.E.H.~Ismail, \emph{Classical and Quantum Orthogonal Polynomials in One Variable}, paperback ed., Cambridge Univ. Press,  2009. 

\bibitem{IsmaK} M.E.H.~Ismail, E.~Koelink, \emph{The $J$-matrix method}, Adv. Appl. Math. 46 (2011), 379-–395.

\bibitem{Koel} E.~Koelink,
\emph{Spectral theory and special functions},  in: \emph{Laredo Lectures on Orthogonal Polynomials and Special Functions}, (eds. R. \'Alvarez-Nodarse, F. Marcell\'an, W. Van Assche),
Nova Sci. Publ., 2004, pp. 45--84.

\bibitem{KoekS} R.~Koekoek, R.F.~Swarttouw, \emph{The Askey-scheme
of hypergeometric orthogonal polynomials and its q-analogue}, online at
\texttt{http://aw.twi.tudelft.nl/\~{}koekoek/askey.html}, Report
98-17, Technical University Delft, 1998. 

\bibitem{Koor-JF} T.H. Koornwinder, \emph{Jacobi functions and  analysis on noncompact semisimple Lie groups}, in: \emph{Special Functions: Group Theoretical Aspects and Applications}, (eds. R.A.~Askey, T.H.~Koornwinder, W.~Schempp), Math. Appl., Reidel,  1984, pp. 1--85.

\bibitem{Koor-LNM} T.H.~Koornwinder, \emph{Special orthogonal polynomial systems mapped onto each other by the Fourier-Jacobi transform}, in: \emph{Orthogonal Polynomials and Applications}, (eds. C. Brezinski, A. Draux, A.P.~Magnus, P.~Maroni, A.~Ronveaux),   LNM 1171, Springer, 1985, pp. 174--183.

\bibitem{Koor1990} T.H.~Koornwinder,  \emph{Jacobi functions as limit cases of $q$q-ultraspherical polynomials},  J. Math. Anal. Appl.  148  (1990), 44--54.

\bibitem{Simo} B.~Simon, 
\emph{The classical moment problem as a self-adjoint finite
difference operator},
Adv. Math. 137 (1998), 82--203.

\bibitem{Vaks} L.L.~Vaksman, \emph{Quantum Bouned Symmetric Domains}, Transl. Math. Monographs 238, Amer. Math. Soc., 2010.

\bibitem{YamaR-PhysRevA} H.A.~Yamani, W.P.~Reinhardt, {$L^2$ discretizations of the continuum: radial kinetic energy, Coulomb Hamiltonian}, Phys. Rev. A {11} (1975) 1144--1156.

\end{thebibliography}
\end{document}